\documentstyle{article}
\catcode`\@=11
\font\cmss=cmss8
\font\cm=cmbx10
\def\alzi{\raise.02truecm\hbox{{\cmss I}}}
\def\ci{\hbox{{\cm 
C}\hskip-.19truecm{\alzi}\hskip.18truecm}}
\font\uyu=cmr10
\font\rmp=cmr12
\let\app=\in
\font\mat=cmr8
\font\matem=cmmi12
\font\mate=cmmi8
\def\min{\hbox{\mate \char60}\hskip1pt}
\def\clopar{\hbox{\matem \char62}}
\def\opar{\hbox{\matem \char60}}
\def\piu{\raise 1pt\hbox{\mat \char43}\hskip1pt}
\def\ugu{\hbox{\mat \char61}\hskip1pt}
\def\men{\hbox{\rmp \char123}\hskip1pt}
\def\minu{\hbox{\mpic \char20}\hskip1pt}
\font\mpic=cmsy8
\def\implica{\hbox{$\,\Rightarrow\,$}}

\newenvironment{equaz}{$$\refstepcounter{teo}}{\eqno{\rm
(\theteo)}$$\global\@ignoretrue}
\newenvironment{numera}{\refstepcounter{teo}}{\hfill{\rm
(\theteo)}\global\@ignoretrue}
\def\eep{\hbox{$\varepsilon$}}
\let\PI=\Pi
\let\ze=\zeta
\def\la{\hbox{$\lambda$}}
\def\De{\hbox{$\Delta$}}
\def\es{\exists\,}
\def\fre{{\longmapsto}}
\def\diverso{\ugu\hskip -7pt / \kern 2pt}
\def\unoenne#1#2{\hbox{$#1_1\vir#1_{#2}$}}
\def\vir{,\ldots ,}
\let\meno=\setminus
\def\ristr#1{\lower4pt\hbox{$|_{#1}$}}
\font\tengt=eufm10 \def\bgo{\hbox{\tengt B}}
\def\modulo#1#2{\setbox1=\hbox{$#2$}\setbox2=\hbox{\cme \char44}
 \raise 3pt \hbox{\hbox{$#1$}
 \hskip-.27truecm\lower.17truecm\copy2\hskip-.2truecm\lower.3truecm\copy1}}
\def\cunoenne#1#2{\hbox{$#1_1\cdots#1_{#2}$}}
\font\cme=line10

\def\Pi{\hbox{{\uyu I}\hskip-.05truecm{\uyu P}}}
\font\ZZZZ=cmssdc10 scaled\magstep1
\def\zeta{\hbox{{\ZZZZ Z}\hspace{-4.3pt}{\ZZZZ 
Z}\kern2pt}}
\def\RE{\rm Re}
\def\IM{\rm Im}
\def\DS{\displaystyle}
\def\F{\hbox{$\cal F\!$}}
\def\mapright#1{\smash{\mathop{\longrightarrow}\limits^{#1}}}
\def\mapdestra{\smash{\mathop{\longrightarrow}}}
\def\FF{\hbox{{\uyu I}\hskip-.05truecm{\uyu F}$\!$}}
\def\puntif{\begin{figure}\centering
\begin{picture}(190,45)(10,0)
 \put(158,17){\makebox(5,10){$P_+$}}\put(160,15){\thicklines\circle{2}}
 \put(126,17){\makebox(5,10){$P_1$}}\put(130,15){\thicklines\circle{2}}
\put(88,21){\makebox(5,10){$P_0$}}\put(90,15){\thicklines\circle{2}}
 \put(58,17){\makebox(5,10){$P_{-1}$}}\put(60,15){\thicklines\circle{2}}
 \put(29,17){\makebox(5,10){$P_-$}}\put(30,15){\thicklines\circle{2}}
\put(-10,15){\line(1,0){190}}\put(104,15){\thicklines\circle{2}}
\put(90,15){\thicklines\circle{10}}\put(98,0){\makebox(5,10){$\gamma_0'$}}
\put(104,15){\thicklines\line(-1,0){10}}
\put(90,20){\thicklines\line(-1,1){4}}\put(90,20){\thicklines\line(-1,-1){4}}
\put(89,8.7){\makebox(3,3){$\scriptstyle >$}}
\end{picture}\\[1em]
\begin{picture}(165,45)(10,0)
 \put(148,17){\makebox(5,10){$P_+$}}\put(150,15){\thicklines\circle{2}}
 \put(116,21){\makebox(5,10){$P_1$}}\put(120,15){\thicklines\circle{2}}
\put(78,17){\makebox(5,10){$P_0$}}\put(80,15){\thicklines\circle{2}}
 \put(48,17){\makebox(5,10){$P_{-1}$}}\put(50,15){\thicklines\circle{2}}
 \put(19,17){\makebox(5,10){$P_-$}}\put(20,15){\thicklines\circle{2}}
\put(0,15){\line(1,0){165}}\put(94,15){\thicklines\circle{2}}
\put(120,15){\thicklines\circle{10}}\put(107,0){\makebox(5,10){$\gamma_1'$}}
\put(95,15){\thicklines\line(1,0){20}}
\put(119,8.7){\makebox(3,3){$\scriptstyle >$}}
  \end{picture}\hspace{1em}\nolinebreak
\begin{picture}(165,45)(10,0)
 \put(148,21){\makebox(5,10){$P_+$}}\put(150,15){\thicklines\circle{2}}
 \put(116,17){\makebox(5,10){$P_1$}}\put(120,15){\thicklines\circle{2}}
\put(78,17){\makebox(5,10){$P_0$}}\put(80,15){\thicklines\circle{2}}
 \put(48,17){\makebox(5,10){$P_{-1}$}}\put(50,15){\thicklines\circle{2}}
 \put(19,17){\makebox(5,10){$P_-$}}\put(20,15){\thicklines\circle{2}}
\put(0,15){\line(1,0){165}}\put(94,15){\thicklines\circle{2}}
\put(150,15){\thicklines\circle{10}}\put(129,0){\makebox(5,10){$\gamma_+'$}}
\put(95,15){\thicklines\line(1,0){20}}\put(125,15){\thicklines\line(1,0){20}}
\put(120,15){\thicklines\oval(10,10)[b]}
\put(149,8.7){\makebox(3,3){$\scriptstyle >$}}
  \end{picture}\\[1em]
\begin{picture}(165,45)(10,0)
 \put(148,17){\makebox(5,10){$P_+$}}\put(150,15){\thicklines\circle{2}}
 \put(116,17){\makebox(5,10){$P_1$}}\put(120,15){\thicklines\circle{2}}
\put(78,17){\makebox(5,10){$P_0$}}\put(80,15){\thicklines\circle{2}}
 \put(48,21){\makebox(5,10){$P_{-1}$}}\put(50,15){\thicklines\circle{2}}
 \put(19,21){\makebox(5,10){$P_-$}}\put(20,15){\thicklines\circle{2}}
\put(0,15){\line(1,0){165}}\put(94,15){\thicklines\circle{2}}
\put(20,15){\thicklines\circle{10}}\put(35,0){\makebox(5,10){$\gamma_-'$}}
\put(95,15){\thicklines\line(-1,0){10}}\put(75,15){\thicklines\line(-1,0){20}}
\put(50,15){\thicklines\oval(10,10)[t]}\put(80,15){\thicklines\oval(10,10)[b]}
\put(45,15){\thicklines\line(-1,0){20}}
\put(19,8.7){\makebox(3,3){$\scriptstyle >$}}
 \end{picture}\hspace{1em}\nolinebreak
\begin{picture}(165,45)(10,0)
 \put(148,17){\makebox(5,10){$P_+$}}\put(150,15){\thicklines\circle{2}}
 \put(116,17){\makebox(5,10){$P_1$}}\put(120,15){\thicklines\circle{2}}
\put(78,17){\makebox(5,10){$P_0$}}\put(80,15){\thicklines\circle{2}}
 \put(48,21){\makebox(5,10){$P_{-1}$}}\put(50,15){\thicklines\circle{2}}
 \put(19,17){\makebox(5,10){$P_-$}}\put(20,15){\thicklines\circle{2}}
\put(0,15){\line(1,0){165}}\put(94,15){\thicklines\circle{2}}
\put(50,15){\thicklines\circle{10}}\put(65,0){\makebox(5,10){$\gamma_{-1}'$}}
\put(95,15){\thicklines\line(-1,0){10}}\put(75,15){\thicklines\line(-1,0){20}}
\put(80,15){\thicklines\oval(10,10)[b]}
\put(49,8.7){\makebox(3,3){$\scriptstyle >$}}
  \end{picture}\caption{\label{figura}A basis of
$\pi_1(\Pi^1\meno\{P_1,P_{-1},P_+,P_-,P_0\})$.}\vspace*{2em}
\end{figure} }
\catcode`\@=12
\begin{document}
\title{The orbifold fundamental group of 
Persson-Noether-Horikawa surfaces.}
\author{Fabrizio Catanese and Sandro Manfredini}
\date{}
\maketitle
\hfill This article is dedicated to the\\
\hbox to .99\textwidth{\hfill memory of Boris 
Moishezon.}
\section{Introduction.}
Among the minimal surfaces of general type, the 
Noether surfaces are
those for which the Noether inequality $K^2 \geq 
2p_g -4$ is an equality
($K^2$ is the self intersection of a canonical 
divisor, $p_g$ is the
dimension of the space of holomorphic 2-forms).\\
These surfaces were described by Noether (\cite{No}) 
and more recently by
Horikawa (\cite{Ho}) who proved that if $8\ |\ K^2$ 
then there are two
distinct deformation types, namely the 
Noether-Horikawa surfaces of
connected type (for short, N-H surfaces of type C), 
and those of non
connected type (for short, of type N). This notation 
refers to the fact
that, the canonical map being a double covering of a 
rational ruled
surface, for type C the branch locus is connected, 
whereas for type N it
is not connected.\\
In particular Horikawa proved that the intersection 
forms are both of the same
parity (in fact, both odd) if and only if $16\ |\  
K^2$.\\
 From M. Freedman's theorem (\cite{Fr}) follows that 
if $16\ |\  K^2$ type N
and type C provide  two orientedly homeomorphic 
compact 4-manifolds.\\
Horikawa posed the question whether type N and type 
C provide
 two orientedly diffeomorphic compact 4-manifolds.\\
It looked like a natural problem to try to see 
whether the two differentiable
structures could be distinguished by means of the 
invariants introduced by
S. Donaldson in \cite{Do}.\\
In the case of type C we have been able (\cite{Ca}) 
to calculate the constant
Do\-nald\-son invariants (corresponding to 
zero-dimensional moduli spaces)
using some singular canonical models of these 
surfaces with very many
singularities, and an approach introduced by P. 
Kronheimer (\cite{Kr}) for
the case of the Kummer surfaces.
The number we obtained, namely $ 2^{2k}$  when $K^2=8k$,
is the leading term  of the Donaldson series (see 
\cite{K-M}), which was
later fully calculated by Fintushel and Stern in the 
case of N-H surfaces
of type C via the technique of rational blow-downs 
(\cite{F-S}).\\
The Donaldson series for N-H surfaces of type N has 
not yet, to our
knowledge, been calculated;
although, after the Seiberg-Witten theory (\cite{W}) 
has been introduced,
and after Pidstrigach and Tyurin (\cite{P-T}) have 
announced the equality
between Kronheimer-Mrowka and Seiberg-Witten 
classes, the two series
should be equal.\\
Our original aim was to extend the application of 
the Kronheimer theory to the
case of N-H surfaces of type N using a very singular 
model constructed by
Ulf Persson (\cite{Per}), describing its orbifold 
fundamental group, its
representations into $SO(3)$, and then trying to see 
which of those have
virtual dimension zero.\\
In this article we consider the singular N-H 
surfaces of type N with maximal
Picard number constructed by Persson, henceforth called
Persson-Noether-Horikawa surfaces (P-N-H for short), 
and we determine
their orbifold fundamental group.\\
This is our main result:\\
{\bf Theorem.} {\em The orbifold fundamental group 
of the P-N-H surfaces is
$$\zeta_4\oplus\zeta_2$$ if $16\ |\ K^2,$
$$\zeta_4\oplus\zeta_4$$\nopagebreak
in the other case where $8\ |\ K^2$  but $16$ does 
not divide
$K^2$.}\\
It follows immediately that we have, for $16\ |\ 
K^2$, only six nontrivial
classes of orbifold $SO(3)$-representations, and a 
result which we do not
prove here is that  we do not get anyone of virtual 
dimension zero.\\
This is not surprising in view of (\cite{P-T}), 
since if Kronheimer's
approach would have worked, we would have had only a 
finite number of
constant Donaldson invariants.\\
On the other hand, the algebro-geometric technique 
of studying canonical
mo\-dels with many rational double points produces 
on the smooth model
configurations of (-2)-projective lines (spheres) 
whose tubular neighborhood
has a unique holomorphic structure and, in 
particular, a unique compatible
$C^{\infty}$ structure. In this way one produces a 
decomposition of the
4-manifold in geometric pieces, one of which is the 
nonsingular part of the
singular canonical model.\\
 From this point of view, the calculation of the 
orbifold fundamental group
leads to a better understanding of the 
differentiable structures of the
smooth model.\\
Since our proof is rather involved technically we 
would like to give a brief
geometrical "explanation" of our result.\\
Persson's construction starts with a plane nodal 
cubic $C$ meeting a conic
$Q$ at only one point $P$. Moreover, $C$ and $Q$ 
have two common tangents
$L_{-1}$ and $L_1$ which meet in a point $O$ 
collinear with $P$ and the node
of $C$.\\
Blowing up $O$ we get a $\Pi^1$-bundle 
$f':\FF_1\mapdestra\Pi^1$ with a
section $\Sigma_{\infty}$, a bisection $Q'$ and a 
3-section $C'$ ($'$
denoting the proper transform under the blow up).\\
A cyclic cover of order $2k\piu2$ branched on 
$L'_{-1}$ and $L'_1$ yields a
new $\Pi^1$-bundle $f'':\FF_{2k+2}\mapdestra\Pi^1$ 
with a section
$\Sigma''_{\infty}$ disjoint from a 3-section $C''$ 
and two sections
$Q_1''$, $Q_2''$ (the inverse image $Q''$ of $Q'$ 
splits into two
components).\\
The curve $B=C''\cup Q_1''\cup 
Q_2''\cup\Sigma''_{\infty}$ has many singular
points, and our canonical model $X_{2k+2}$ is the 
double cover of
$\FF_{2k+2}$ branched on $B$. By construction 
$X_{2k+2}$ has a genus 2
fibration onto $\Pi^1$, whence the orbifold 
fundamental group
$\pi_1(X_{2k+2}^{\#})$, $X_{2k+2}^{\#}$ being the 
nonsingular part of
$X_{2k+2}$, is a quotient of $\pi_1(F)$, where $F$ 
is a fixed genus 2
fibre.\\
$F$ being a double cover of $\Pi^1$ branched in six 
points
$P_0\ugu\Pi^1\cap\Sigma''_{\infty}$,
$P_1\ugu\Pi^1\cap Q''_1$, $P_2\ugu\Pi^1\cap Q''_2$,
$\{P_3,P_4,P_5\}=\Pi^1\cap C''$, $\pi_1(F)$ is the 
subgroup of a free
product $\F_5(2)$ of five copies of $\zeta_2$, given 
by words of even
length.\\
$\F_5(2)$ is generated by elements 
$\unoenne{\eep}{6}$ such that
$\cunoenne{\eep}{6}\ugu1$ ($\eep_i$ corresponds to a 
loop in $\Pi^1$ around
the point $P_{i-1}$).\\
The first main point (we must be rather vague here, 
else we must give the
full proof) is that, since curve $C''$ is 
irreducible, when the fibre $F$
moves around, $\eep_4,\eep_5,\eep_6$ become 
identified.\\
Thus we only have $\unoenne{\eep}{4}$ with 
$\cunoenne{\eep}{4}\ugu1$, and
therefore we have "proved" that our group is abelian 
, being a quotient of
the fundamental group $\Gamma$ of a curve of genus 1 
obtained as the double
cover of $\Pi^1$ branched in four points. More 
precisely, $\Gamma$ is an
abelian group with generators $\eep_1\eep_2$, 
$\eep_1\eep_3$.\\
We must still take into account the fact that, when 
the fibre $F$ moves
towards a singular point (corresponding to points of 
intersection $C''\cap
Q_1''$, $C''\cap Q_2'', Q_1''\cap Q_2''$), further 
relations are introduced.
These relations are hard to control globally but if 
we look locally
around these points of intersection, and accordingly 
take a new basis
$\unoenne{\eep'}{4}$, the situation becomes simpler.\\
In fact, the local equation of the double cover is 
$z^2\ugu y^2\men x^{2c}$,
where $c\ugu6$ or $c\ugu k+1$, and $x$ is the 
pullback of a local coordinate
on $\Pi^1$, so that the corresponding local braid 
yields the relation
$(\eep'_j\eep'_i)^c=(\eep'_i\eep'_j)^c$. In turn, 
using $(\eep'_i)^2\ugu1$,
we obtain the relation $(\eep'_j\eep'_i)^{2c}\ugu1$.\\
That's how one shows that the two generators of the 
abelian group have
period 2 or 4.\\
The paper is organized as follows:\\
In section two we take up Persson's construction 
using explicit
equations showing that the surface is defined over a 
real quadratic field.\\
In the third section we describe the five steps 
leading to a presentation
of our fundamental group in terms of the braid 
monodromy of the plane
curve $D=C\cup Q$.\\
Finally, in section four we apply combinatorial 
group theory arguments in
order to give the main result concerning the 
orbifold fundamental group.\\
Acknowledgements : Both authors acknowledge support 
from the AGE Project
H.C.M. contract ERBCHRXCT 940557 and from 40\% 
M.U.R.S.T..\\
The first author would like to express his gratitude 
to the Max-Planck
Institut in Bonn where this research was initiated 
(in 1993), and to the
Accademia dei Lincei where he is currently 
Professore Distaccato.
\section{Persson's configuration.}
In this section we will provide explicit equations 
for the configuration
constructed by Ulf Persson in \cite{Per}.\\This is 
the configuration formed
by a smooth conic $Q$ and a nodal cubic $C$ 
intersecting in only one point
$P$ which is smooth for $C$. Moreover $Q$ and $C$ 
have two common tangents
$L_{1}$ and $L_{-1}$ meeting in a point $O$ lying on 
the line joining $P$
and the node of $C$.\\ Let $Q\subset\ci\Pi^2$
be the conic $\{(x,y,z)\app\ci\Pi^2\,|\,x^2\piu 
2zy\piu z^2\ugu 0\}$.\\
Since $$x^2\piu 2zy\piu z^2\ugu(x\piu z)^2\piu 
2z(y\men x)\ugu (x\men
z)^2\piu 2z(y\piu x)$$
$Q$ is tangent to the lines $L_1=\{x\men y\ugu 0\}$ and
$L_{-1}=\{x\piu y\ugu 0\}$.\\
The tangency points are:
$$ x\men y\ugu x\piu z\ugu 0 \implica (1,1,\men1)$$
$$x\piu y\ugu x\men z\ugu 0\implica (1,\men1,1).$$
Note that $Q$ is also tangent to the line $z\ugu 0$ 
at the point
$(0,1,0)=P$.\\
We want to find an irreducible nodal cubic $C$ such 
that $C\cdot Q=6P$ and
such that $C$ is tangent to the lines $x\ugu\pm y$ 
in points different from
those of $Q$.\\
Let $C$ be a cubic s.t. $P\app C$ and $C\cdot Q=6P$. 
Note that if $C$ were
reducible, then the previous condition would imply 
that $z\ugu0$ is a
component of $C$.\\
We then have ${\rm div}(C)={\rm div}(z^3)\ ({\rm 
mod}Q)$, so $C=z^3\piu QL$
with $L$ a linear form, and thus
$$C=z^3\piu (x^2\piu 2zy\piu z^2)(ax\piu by\piu cz).$$
Since we want $C$ to be tangent to the two lines 
$L_{1}$ and $L_{-1}$
we obtain that the following homogeneous polynomials 
in $(x,z)$
\begin{equaz}z^3\piu (x\piu z)^2((a\piu
b)(x\piu z)\piu z(c\men a\men b))\label{1}\end{equaz}
\begin{equaz}z^3\piu (x\men
z)^2((a\men b)(x\men z)\piu z(c\piu a\men 
b))\label{21}\end{equaz} must have
a double root.\\
Set $\ze\ugu (\frac{z}{x+z})$ and $\hat\ze\ugu 
(\frac{z}{x-z})$ and rewrite
\ref{1}, \ref{21} as:
$$\ze^3\piu \ze (c\men a\men b)\piu (a\piu b)=0\ \ \ 
\ \hat\ze^3\piu\hat\ze
(c\piu a\men b)\piu (a\men b)=0.$$
We recall that if $\ze$ is a double root of $z^3\piu 
pz\piu q\ugu 0$ then
$$3\ze^2+p=0{ \rm\ \ whence\ \ } \frac23\ze p\piu 
q\ugu 0$$ and this implies
that $$\ze\ugu\men \frac32\frac pq{\rm\ \ thus\ \ 
}27q^2\piu 4p^3\ugu 0.$$
Therefore we have a double root of \ref{1} if and 
only if $$\es A\ :\
\ze\ugu
A,\ q\ugu 2A^3,\ p\ugu\men 3A^2,\ {\rm i.e.}\ 
\left\{\begin{array}{l} a\piu
b\ugu 2A^3 \\ c\men (a\piu b)\ugu\men 
3A^2.\end{array}\right.$$
Similarly if we set $\hat\ze\ugu\men B$
we have $$\left\{\begin{array}{l} a\men b\ugu\men 2B^3\\
c\men (b\men a)\ugu\men3B^2\end{array}\right.$$ and so
$$\left\{\begin{array}{l}a\ugu A^3\men B^3 \\ b\ugu 
A^3\piu B^3 \\
c\ugu 2A^3\men 3A^2\ugu 2B^3\men 3B^2. 
\end{array}\right.$$
Then $A$ and $B$ must satisfy $2(A^3\men B^3)\ugu 
3(A^2\men B^2)$.\\
Recall that (we make no distinction between a curve 
and its equation)
$$C=z^3\piu (x^2\piu 2zy\piu z^2)((A^3\men
B^3)x\piu (A^3\piu B^3)y\piu (2A^3\men 3A^2)z)$$
while $x\men y\ugu 0$ is tangent to $C$ at the point 
where
$$\ze\ugu\frac{z}{x+z}\ugu A.$$  Therefore the 
tangency point is
$(1\men A,1\men A,A)$.\\ Similarly $x\piu y\ugu 0$ 
is tangent to $C$ at the
point $(B\men 1,1\men B,B)$.\\
Let us now search for a cubic $C$ with a singular 
point on the line
$x\ugu0$, as in Persson's construction.\\
Since $\frac{\partial C}{\partial x}$ on the line 
$x\ugu0$ equals $aQ$ and
the singular point is different from $P$ it follows that
$a\ugu0$. Whence $A^3\men B^3\ugu A^2\men B^2\ugu 0$ 
and so $A\ugu B$.\\
If $A\ugu B$ then $C$ contains only the monomial 
$x^2$ as a polynomial in
$x$, so the involution $x\fre \men x$
leaves the curve $C$ invariant. From this we deduce 
that a singular point of
$C$ must have its $x$ coordinate equal to $0$ and 
$C$ has then a
singularity on the line $x\ugu 0$ if and only if
$$z^3\piu (z^2\piu 2zy)(2A^3y\piu (2A^3\men 3A^2)z)\ 
\ {\rm has\ a\ double\
root.}$$ Remembering that it can't be $A\ugu B\ugu 
0$, the double root
cannot be $z\ugu 0$ and we can write the above as 
$$z(z^2\piu (z\piu
2y)(2A^3y\piu (2A^3\men 3A^2)z)).$$ So we must check 
that
$$z^2(1\piu 2A^3\men 3A^2)\piu 2zy(A^3\piu 2A^3\men 
3A^2)\piu 4A^3y^2=$$
$$=z^2(1\piu 2A^3\men 3A^2)\piu 2zy3A^2(A\men 1)\piu 
4A^3y^2$$
has a double root.\\ This is the case when
$$9A^4(A\men 1)^2\ugu 4A^3(1\piu 2A^3\men 3A^2)\ \ 
{\rm\ i.e.}$$
$$9A^6\men 18A^5\piu 9A^4\ugu 4A^3\piu 8A^6\men 12A^5.$$
Upon dividing by $A^3\diverso 0$ we get
$$A^3\men 6A^2\piu 9A\men 4\ugu 0.$$
Observe that $1$ is a root of this equation, but if 
$A\ugu 1$ then the
singular point is $(0,0,1)$ and coincides with the 
point of tangency of
$x\piu y\ugu 0$ so this root has to be discarded. Since
$$A^3\men 6A^2\piu 9A\men 4\ugu (A\men 1)(A^2\men 
5A\piu 4)\ugu (A\men
1)^2(A\men 4)$$
the other possible root is then $A\ugu 4$, and in 
this case we have $B\ugu
A\ugu 4$, $a\ugu 0$, $b\ugu 8\cdot4^2$, $c\ugu 
5\cdot4^2$.\\
Then
$$C=z^3\piu 4^2(x^2\piu 2yz\piu z^2)(8y\piu 5z)$$
The tangency points are $(\men 3,\men 3,4)$ and 
$(3,\men 3,4)$, while for
the singular point we have $x\ugu 0$ and a double 
root of
$$z^2\piu 4^2(2y\piu z)(8y\piu 5z)\ugu 0\iff 
81z^2\piu 4^218zy\piu
4^4y\ugu 0\iff 9z\piu 4^2y\ugu 0$$
so the singular point is $(0,9,\men 16)$.\\
With this choice of $A$ and $B$, $C$ is irreducible 
(since $z\ugu 0$ is
not a component of $C$).\\
We want to find the lines through $(0,0,1)$ and 
tangent to $C$.\\
Let $A\ugu B\ugu\la$ and consider more generally the 
1-parameter family
of curves: $$C_{\la}=z^3\piu (x^2\piu 2yz\piu 
z^2)(2\la^3y\piu (2\la^3\men
3\la^2)z)\ugu 0.$$
The tangency points on the two fixed lines $x\piu 
y=0$, $x\men y=0$ are, as
we know, $(1\men\la,1\men\la,\la)$ and $(\la\men 
1,1\men\la,\la)$.\\
Rewriting the last equation in powers of $z$ we obtain:
$$z^3(1\piu 2\la^3\men 3\la^2)\piu 
z^26y\la^2(\la\men 1)\piu z\la^2(4\la
y^2\piu (2\la\men 3)x^2)\piu 2\la^3x^2y\ugu 0.$$
Since we know what happens for $\la\ugu 0$, we can 
divide by $\la^3$, set
$w\ugu \frac z{\lambda}$ and obtain:
$$w^3(1\piu 2\la^3\men 3\la^2)\piu w^26y\la(\la\men 
1)\piu w(4\la
y^2\piu (2\la\men 3)x^2)\piu 2x^2y\ugu 0.$$
We let now $\De$ be the discriminant of $C_{\la}$ 
with respect to the
variable $w$, and using a
standard formula for $\De$, we find a degree 6 
equation in $x$ and $y$ which
is divisible by $x^2(x^2\men y^2)$.\\
Remembering that the discriminant of $a_0x^3\piu 
a_1x^2\piu a_2x\piu a_3$
is: $$\De\ugu a_1^2a_2^2\men 4a_0a_2^3\men
4a_1^3a_3\men 27a_0^2a_3^2\piu 18a_0a_1a_2a_3$$
and applying this formula for simplicity when 
$\la\ugu 4$, we obtain:
$$y^22^63^4(16y^2\piu 5x^2)^2\men 2^23^4(16y^2\piu 
5x^2)^3\men$$
$$\men2^{12}3^6x^2y^4 \men 2^23^{11}x^4y^2\piu 
2^53^8(16y^2\piu
5x^2)x^2y^2$$
and factoring this binary form we get:
$$x^2(x^2\men y^2)2^23^4(2^7y^2\men 5^3x^2).$$
So we have that the tangent lines to $C$ passing 
through $(0,0,1)$ are
$x\ugu\pm y$, $x\ugu\pm\sqrt{\frac{128}{125}}y$ 
while $x\ugu 0$ passes
through the node of $C$. We denote by $L_0$ the line 
$x\ugu0$ and by
$L_+,L_-$ the two lines $x\ugu\sqrt{\frac{128}{125}}y$,
$x\ugu-\sqrt{\frac{128}{125}}y$ respectively.\\
In order to find the tangency point on the lines 
$L_+,L_-$
we by symmetry may restrict to the line $L_+$.\\
Writing $x\ugu 2^3\sqrt2\,a$,
$y\ugu 5\sqrt5\,a$ we have that
\begin{equaz}\label{pol}z^3\piu 2^4(2^7a^2\piu 
10\sqrt5\,az\piu
z^2)(40\sqrt5\,a\piu
5z)\ugu 0\end{equaz} has a double root. Since for 
its derivative we have
$$3z^2\piu 2^4(10\sqrt5\,a\piu 2z)(40\sqrt5\,a\piu 
5z)\piu 2^45(2^7a^2\piu
10\sqrt5\,az\piu z^2)\ugu 0$$
$$(15\piu {3\over16})z^2\piu 180\sqrt5\,az\piu 
2640a^2\ugu0$$
$${a\over z}\ugu{-90\sqrt5\pm\sqrt{90^25-
2640(15+{3\over16})}\over2640}\ugu\sqrt5{-30\pm3\over880}.$$
Thus $\frac{y}{z}={-25(30\pm3)\over880}$,
$\frac{x}{z}={-8\sqrt{10}(30\pm3)\over880}$ and the 
point of tangency is
one of the points $(\men33\cdot8\sqrt{10},\men 
25\cdot 33,880)$,
$(\men27\cdot8\sqrt{10},\men 25\cdot 27,880)$.\\
Upon substituting these values in the polynomial 
\ref{pol} we find that the
correct choice is $(\men24\sqrt{10},\men 75,80)$.\\ 
By symmetry the point
$(24\sqrt{10},\men 75,80)$  is the tangency point of 
the line $L_-$.\\
Let us write
$$C=4^2(8y\piu 5z)x^2\piu z(16y\piu 9z)^2\ugu 0$$
and let us set $u\ugu 16y\piu 9z$. We have:
$$C= zu^2\piu 8x^2(u\piu z)\ugu 0$$
In these coordinates the singular point of $C$ is 
$(0,0,1)$, so the
tangents at the singular point are given by:
$$8x^2\piu u^2\ugu 0$$
whence they are complex and we have an isolated point.\\
In order to draw $C$, let's compute its flexes. 
Using the
coordinates $x$, $u$, and $z$ the Hessian matrix is:
$$\left(\begin{array}{ccc}16(u\piu
z)&16x&16x\\16x&2z&2u\\16x&2u&0\end{array} \right)$$
The Hessian curve is then given by the determinant of
$$\left(\begin{array}{ccc}(u\piu z)&0&x\\0&z\men 
2u&u\\8x&u&0\end{array}
\right)$$
which equals
$$\men (u\piu z)u^2\men 8x^2(z\men 2u)\ugu0.$$\\
Eliminating $8x^2$ from the two equations we get
$$(u\piu z)^2u^2\men zu^2(z\men 2u)\ugu0$$
so either $u=0$, and this implies either $x\ugu 0$ 
(the singular point) or
$z\ugu 0$ that gives the point $(1,0,0)$, or
$$(u\piu z)^2\men z(z\men 2u)\ugu u^2\piu 4uz\ugu 0$$
that gives ($u\diverso0$) $u\ugu \men 4z$, that is 
$z\ugu\men 1$, $u\ugu 4$,
$y\ugu{13\over16}$, $x\ugu\pm\sqrt{2/3}$.\\
For these points ${x\over 
y}=\pm\sqrt{2/3}{16\over13}$.\pagebreak
\section{Fundamental groups.}
In this section we are going to describe the five 
steps leading to the
determination of the orbifold fundamental group of 
the Persson's surfaces.\\
{\bf Step 1.}\\
Let $\FF_1$ be the blow up of $\Pi^2$ at the point 
$(0,0,1)$ and let
$\Sigma_{\infty}$ be the exceptional divisor.\\
We consider the fibre bundle 
$\FF_1\mapright{f'}\Pi^1$ and its
restriction $f$ $$\FF_1\meno( C\cup 
Q\cup\Sigma_{\infty}\cup L_1\cup L_{-1}
\cup L_+\cup L_-\cup L_0)=\tilde{\FF}_1$$
$$f\big\downarrow$$
$$\Pi^1\meno\{P_1,P_{-1},P_+,P_-,P_0\}=\Pi^1\meno\{5\ {\rm
pts.}\}.$$ $f$ is again a fibre bundle and we have a 
corresponding homotopy
exact sequence of fundamental groups $$1\mapdestra
\F_5\mapdestra\tilde{\PI}\mapdestra\F_4\mapdestra1$$
\puntif
where $\F_k$ denotes the free group with $k$ 
generators and
$\tilde\PI=\pi_1(\tilde{\FF}_1)$.\\
Here we choose a small positive real number 
$\eep\mag0$ and $x\ugu\eep$,
$y\ugu1$ as base point on $\Pi^1\meno\{5\ {\rm 
pts.}\}$ and $x\ugu\eep$,
$y\ugu1$, $z\ugu\men4\sqrt{-1}$ as base point on 
$\FF_1$.\\
We let $\unoenne{\delta}{5}$ be a natural geometric 
basis of the free
group $$\F_5=\pi_1(f^{-1}({\rm base\ 
pt.}))=\pi_1(L_{\eep}\meno(C\cup
Q\cup\Sigma_{\infty}))$$ where the five points 
$L_{\eep}\cap C$,
$L_{\eep}\cap
Q$ are ordered by lexicografic order on $\RE(\frac 
zy)$, $\IM(\frac zy)$.\\
$\F_4$ is generated by the five geometric paths 
$\gamma_i'$ around the five
critical values described in figure \ref{figura} and 
whose product is the
identity.\\
For these elements we choose lifts to $\FF_1$ using 
a $C^{\infty}$ section
of a tubular neighborhood of $\Sigma_{\infty}$ 
meeting $\Sigma_{\infty}$
just in the point $\infty$ $(y\ugu0)$ with 
intersection number equal to
$-1$.\\
Therefore such lifts give paths $\gamma_i$ such that
$$\prod\gamma_i=\prod\delta_i$$ and more specifically
$$\gamma_+\gamma_1\gamma_0\gamma_-\gamma_{-1}=\delta_1\cdots\delta_5=
\gamma_{-1}\gamma_+\gamma_1\gamma_0\gamma_-.$$
We have that, indeed, $\tilde{\PI}$ occurs as a 
semidirect
product described by the relations
$$\gamma_j^{-1}\delta_i\gamma_j=(\delta_i)\beta_j$$
where the $\beta_j$'s are suitable braids in
$$\bgo_5
  =\opar\sigma_1,\ldots,\sigma_{4}|\:\:\:\,
  \sigma_i\sigma_j \ugu\sigma_j\sigma_i\ \;\forall\: 
1\minu
  i\min j\piu1\minu 5
$$
\hspace*{15.4em}\immediate\vspace*{-1ex}
$\sigma_i\sigma_{i+1}\sigma_i\ugu 
\sigma_{i+1}\sigma_i\sigma_{i+1}\
\forall\: 1\minu i\min 4\ \clopar$
\vspace{\baselineskip}\\
the braid group on 5 strings which acts on the right 
on the free group
$\F_5$
by the formulae 
\begin{eqnarray*}(\delta_h)\sigma_k&=&\delta_h\ \ \ 
\ {\rm
if}\ h\diverso k,k\piu1\\
(\delta_k)\sigma_k&=&\delta_{k+1}\\
(\delta_{k+1})\sigma_k&=&\delta_{k+1}^{-1}\delta_k\delta_{k+1}.
\end{eqnarray*}
The braids $\beta_j$ are constructed by following 
the motion of the five
points of the intersection of $f'^{-1}(P)$ with 
$C\cup Q$ while $P$ goes
along $\gamma_j'$.\\
With our choice of the $\gamma'$'s we have, as the 
reader can easily verify,
\begin{eqnarray*}
\beta_0&=&\sigma_4^{12}\sigma_2^2\\
\beta_1&=&\sigma_1^{-1}\sigma_2\sigma_3\sigma_1\sigma_2^{-1}\sigma_1\\
\beta_{-1}&=&\sigma_4^{-6}\sigma_2^{-1}\beta_1\sigma_2\sigma_4^6\\
\beta_+&=&\sigma_1^{-2}\sigma_2\sigma_3\sigma_4\sigma_3^{-1}\sigma_2^{-1}
\sigma_1^2\\
\beta_-&=&\sigma_4^{-6}\sigma_2^{-1}\beta_+\sigma_2\sigma_4^6.
\end{eqnarray*}
{\bf Step 2.}\\
By taking $\sqrt{\frac{x-y}{x+y}}$ we have a new 
fibre bundle
$\FF_2\mapright{g'}\Pi^1$ obtained by base change. 
Under this base change
the inverse image $Q'$ of the conic $Q$ splits into 
two sections of $g'$
which we will denote by $Q_1'$ and $Q_2'$. Again, by 
restriction we have a
fibre bundle $g$
$$\hat{\FF}_2=\FF_2\meno (C'\cup Q_1'\cup Q_2'\cup 
\Sigma_{\infty}'
\cup \{8\ {\rm fibres}\})\mapright{g}\Pi^1\meno \{8\ 
{\rm
pts.}\}.$$ Correspondingly we get an exact sequence
$$1\mapdestra\F_5=\opar\unoenne{\delta}{5}\clopar\mapdestra\hat{\PI}
\mapdestra\F_7=
\opar\gamma_0,\gamma_-,\gamma_+,\bar{\gamma}_0,\bar{\gamma}_-,
\bar{\gamma}_+,\gamma_1^2\clopar\mapdestra1$$
where 
$\bar\gamma_i\ugu\gamma_i^{\gamma_1}\ugu\gamma_1\gamma_i\gamma_1^{-1}$
and $\hat\PI=\pi_1(\hat{\FF}_2)$.\\
The fact that $\F_7$ has seven generators as above 
follows since the double
cover of $\Pi^1\meno\{5\ {\rm pts.}\}$
corresponds to the homomorphism 
$\F_4\rightarrow\zeta_{\!2}$ sending
$\gamma_1',\gamma_{-1}'\mapsto\bar1$, and 
$\gamma_0',\gamma_+',
\gamma_-'\mapsto\bar0$.\\
If we want to keep track of the eight critical 
values, we can also use
$(\gamma_{-1}^2)^{\gamma_1}$ as a generator. In fact
$$(\delta_1\cdots\delta_5)^2=(\gamma_+\gamma_1\gamma_0\gamma_-\gamma_{-1})
(\gamma_{-1}\gamma_+\gamma_1\gamma_0\gamma_-)$$
thus
$$\gamma_+\gamma_0^{\gamma_1}\gamma_-^{\gamma_1}(\gamma_{-1}^2)^{\gamma_1}
\gamma_+^{\gamma_1}\gamma_1^2\gamma_0\gamma_-=(\delta_1\cdots\delta_5)^2.$$
The geometric meaning of the above formula is 
related to the fact that
$(\Sigma_{\infty}')^2=-2,$ and more precisely to the 
fact that the new
generators of $\F_7$ lie in a $C^{\infty}$ section 
meeting
$\Sigma_{\infty}'$ in one point with intersection 
number $(-2)$, and not
meeting the other curves.\\
A presentation of $\hat{\PI}$ is thus given by
$$\opar\unoenne{\delta}{5},\gamma_0,\gamma_-,\gamma_+,
\bar{\gamma}_0,\bar{\gamma}_-,\bar{\gamma}_+,\Gamma\ugu\gamma_1^{-2}\ |\
\gamma_0^{-1}\delta_i\gamma_0=(\delta_i)\beta_0\ \ \ 
\ \ \ \ \ \ \ $$
$$\hspace*{16.5em}\vdots$$
$$\hspace*{18.5em}\bar{\gamma}_0^{-1}\delta_i\bar{\gamma}_0=(\delta_i)
\beta_1\beta_0\beta_1^{-1}$$ $$\hspace*{16.5em}\vdots$$
$$\hspace*{17em}\Gamma\delta_i\Gamma^{-1}=(\delta_i)\beta_1^2\clopar$$
{\bf Step 3.}\\
The fundamental group
$$\PI'=\pi_1(\FF_2\meno (C'\cup Q_1'\cup Q_2'\cup 
\Sigma_{\infty}'
\cup L_1'\cup L_{-1}'))$$ is a quotient of $\hat\PI$.
The presentation of $\PI'$ is readily accomplished 
simply by introducing in
the above presentation the further relations
$$\gamma_0=\gamma_-=\gamma_+=\bar{\gamma}_0=\bar{\gamma}_-=
\bar{\gamma}_+=1.$$
Then $\PI'$ is presented as
$$\opar\unoenne{\delta}{5},\Gamma\ |\
\delta_i=(\delta_i)\beta_0\ \ \
\delta_i=(\delta_i)\beta_-\ \ \
\delta_i=(\delta_i)\beta_+$$
$$\hspace*{7em}\delta_i=(\delta_i)\beta_1\beta_0\beta_1^{-1}\ \ \
\delta_i=(\delta_i)\beta_1\beta_-\beta_1^{-1}$$
$$\hspace*{8em}\delta_i=(\delta_i)\beta_1\beta_+\beta_1^{-1}\ \ \
\Gamma\delta_i\Gamma^{-1}=(\delta_i)\beta_1^2\clopar$$
{\bf Remark:} with the new relations we get, setting
$\Gamma_{-1}=(\gamma_{-1}^2)^{\gamma_1},$
$$\Gamma_{-1}\Gamma=(\delta_1\cdots\delta_5)^2$$
{\bf Step 4'}.\\
We denote by $X_2^{\#}$ the non singular
part of the double cover $X_2$ of $\FF_2$ (branched 
over $C'\cup Q_1'\cup
Q_2' \cup\Sigma_{\infty}'$) and by $Z_2^{\#}$ the 
complement
in $X_2^{\#}$ of $L_1'',L_{-1}''$, the respective 
inverse images of
$L_1',L_{-1}'$.\\ We finally let $Y_2^{\#}$ be the 
double cover of
$\FF_2\meno (C'\cup Q_1'\cup Q_2'\cup 
\Sigma_{\infty}' \cup L_1'\cup
L_{-1}')$.\\ Thus $Y_2^{\#}\subset Z_2^{\#}\subset 
X_2^{\#}$.\\
Clearly $\pi_1(Y_2^{\#})=\ker(\PI'\mapdestra\zeta_2)$,
where $\delta_i\fre\bar1$ and $\Gamma\fre\bar0$, is 
generated by
$\Gamma$, $\sigma=\delta_1\Gamma\delta_1^{-1}$, 
$A_i=\delta_1\delta_i$
$(i\ugu1\vir 5)$ and $B_j=\delta_j\delta_1^{-1}$ 
$(j\ugu2\vir 5)$.\\
To find the relations we apply the 
Reidemeister-Shreier rewriting process
to the relations $R_{\alpha}$ of $\PI'$ and to the 
relations
$\delta_1R_{\alpha}\delta_1^{-1}$.\\
{\bf Step 4''}.\\
Clearly, $\pi_1(Y_2^{\#})$ maps onto 
$\pi_1(Z_2^{\#})$ surjectively with
kernel normally generated by $\delta_1^2$, 
$\delta_i^2=B_iA_i$ $(i\ugu2\vir
5)$ and $(\delta_1\cdots\delta_5)^2$, thus
$\pi_1(Z_2^{\#})$ is generated by $A_2\vir 
A_5,\Gamma$ and has for
relations the relations coming from the rewriting of 
$R_{\alpha}$,
$\delta_1R_{\alpha}\delta_1^{-1}$, and the rewriting of
$(\delta_1\cdots\delta_5)^2=1\!$, i.e. 
$A_2A_3^{-1}A_4A_5^{-1}
A_2^{-1}A_3A_4^{-1}A_5=1$.\\
{\bf Remark:} This relation says that the four 
generators $A_2\vir A_5$
are the generators of $\pi_1({\rm fibre})=\pi_1({\rm 
genus\ 2\ curve})$.\\
{\bf Step 5.}\\ Let $m\ugu k\piu1$ and
consider $X_{2m}^{\#}$, the non singular part of the 
m-fold cyclic cover of
$X_2$ totally branched over $L_1''$ and $L_2''$.\\
To find a presentation of $X_{2m}^{\#}$ we first 
need a presentation
of the kernel of the map 
$\pi_1(Z_{2}^{\#})\mapdestra \zeta_m$ such that
$A_i\fre\bar0$ and $\Gamma,\sigma\fre\bar 1$, and 
then we add the relations
$\Gamma^m=\sigma^m=1$.\\
Applying the Reidemeister-Shreier method, we find 
that the kernel
is generated by $\Gamma^m$, $\Gamma^iA_j\Gamma^{-i}$ 
for $i\ugu1\vir m\men1$
and $j\ugu2\vir 5$, by $\Gamma^i\sigma\Gamma^{-i-1}$ 
for $i\ugu1\vir m\men2$
and $\Gamma^{m-1}\sigma$; it has for relations the 
rewriting in the
new generators of the relations $R_{\alpha}''$ of 
$\pi_1(Z_2^{\#})$
and the rewriting of 
$\Gamma^iR_{\alpha}''\Gamma^{-i}$ for
$i\ugu1\vir m\men1$.
\section{Calculations.}
{\bf Step 3.}\\
We have \begin{eqnarray*}
\beta_0&=&\sigma_4^{12}\sigma_2^2\\
\beta_1&=&\sigma_1^{-1}\sigma_2\sigma_3\sigma_1\sigma_2^{-1}\sigma_1\\
\beta_{-1}&=&\sigma_4^{-6}\sigma_2^{-1}\beta_1\sigma_2\sigma_4^6\\
\beta_+&=&\sigma_1^{-2}\sigma_2\sigma_3\sigma_4\sigma_3^{-1}\sigma_2^{-1}
\sigma_1^2\\
\beta_-&=&\sigma_4^{-6}\sigma_2^{-1}\beta_+\sigma_2\sigma_4^6
\end{eqnarray*}
The relations $\delta_i\ugu(\delta_i)\beta_0$ are 
equivalent to the two
relations \begin{equaz}\label{fofififo}
(\delta_4\delta_5)^6=(\delta_5\delta_4)^6\end{equaz}
\begin{equaz}\label{tttt}
\delta_2\delta_3=\delta_3\delta_2.\end{equaz}
The relations $\delta_i\ugu(\delta_i)\beta_+$ amount to
\begin{equaz}\label{totot}
\delta_5=\delta_2^{-1}\delta_1^{-1}\delta_2\delta_1\delta_2.\end{equaz}
In fact, here and in the sequel, we use the 
following argument:
$\beta_+$ is a conjugate $\sigma\sigma_4\sigma^{-1}$ 
of the braid $\sigma_4$
and the braid $\sigma_4$ yields the relation 
$\delta_4\ugu\delta_5$.
Therefore, if we set $\delta_4'\ugu 
(\delta_4)\sigma^{-1}$, $\delta_5'\ugu
(\delta_5)\sigma^{-1}$, we get the relation 
$\delta_4'\ugu\delta_5'$.
By our particular choice of $\sigma$ \begin{eqnarray*}
\delta_5'&\ugu&\delta_5\\
\delta_4'&\ugu&(\delta_4)\sigma_3^{-1}\sigma_2^{-1}\sigma_1^2\\
&\ugu&(\delta_3)\sigma_2^{-1}\sigma_1^2\\
&\ugu&(\delta_2)\sigma_1^2\\&\ugu&(\delta_2^{-1}\delta_1\delta_2)\sigma_1\\
&\ugu&\delta_2^{-1}\delta_1^{-1}\delta_2\delta_1\delta_2\end{eqnarray*}
Similarly, the relations 
$\delta_i\ugu(\delta_i)\beta_-$ are equivalent to
the relation \begin{equaz}\label{thothoth}
\delta_3^{-1}\delta_1^{-1}\delta_3\delta_1\delta_3=
(\delta_4\delta_5)^{-3}\delta_5(\delta_4\delta_5)^3.
\end{equaz}
We write down, for convenience of the reader, the 
action of the braid
$\beta_1^{-1}$, since the new relations $\delta_i
\ugu(\delta_i)\beta_1\beta_j\beta_1^{-1}$ will be 
obtained from the
relations equivalent to 
$\delta_i\ugu(\delta_i)\beta_j$ simply by applying
the automorphism $\beta_1^{-1}$.\begin{eqnarray*}
(\delta_1)\beta_1^{-1}&\ugu&
\delta_1\delta_2\delta_3\delta_2^{-1}\delta_1^{-1}
\delta_2^{-1}\delta_1\delta_2\delta_4\delta_2^{-1}\delta_1^{-1}
\delta_2\delta_1\delta_2\delta_3^{-1}\delta_2^{-1}\delta_1^{-1}\\
(\delta_2)\beta_1^{-1}&\ugu&
\delta_1\delta_2\delta_3\delta_2^{-1}\delta_1^{-1}\\
(\delta_3)\beta_1^{-1}&\ugu&\delta_2^{-1}\delta_1\delta_2\delta_4^{-1}
\delta_2^{-1}\delta_1^{-1}\delta_2\delta_1\delta_2\delta_4\delta_2^{-1}
\delta_1^{-1}\delta_2\\
(\delta_4)\beta_1^{-1}&\ugu&\delta_2^{-1}\delta_1\delta_2\\
(\delta_5)\beta_1^{-1}&\ugu&\delta_5
\end{eqnarray*}
Thus, the relations 
$\delta_i\ugu(\delta_i)\beta_1\beta_0\beta_1^{-1}$
are equivalent to the relations 
\begin{equaz}\label{otwtwo}
(\delta_1\delta_2)^6=(\delta_2\delta_1)^6\end{equaz}
\begin{equaz}\label{ftffff}
\delta_5\delta_3\delta_5^{-1}\delta_4^{-1}\delta_5\delta_4=
\delta_4^{-1}\delta_5\delta_4\delta_5\delta_3\delta_5^{-1}
\end{equaz} where we have used \ref{totot}.\\
The relations 
$\delta_i\ugu(\delta_i)\beta_1\beta_+\beta_1^{-1}$ are
equivalent to the relation 
\begin{equaz}\label{ftotfftfftot}
\delta_5=\delta_2^{-1}\delta_1\delta_2\delta_4^{-1}\delta_5
\delta_3\delta_5^{-1}\delta_4\delta_2^{-1}\delta_1^{-1}\delta_2\end{equaz}
where we have used \ref{totot}.\\
The relations 
$\delta_i\ugu(\delta_i)\beta_1\beta_-\beta_1^{-1}$ 
pop up to
the relation \begin{equaz}\label{thototh}
\delta_3^{-1}\delta_1\delta_2\delta_1^{-1}\delta_3=
(\delta_4\delta_5)^{-2}\delta_5(\delta_4\delta_5)^2.
\end{equaz} In fact, using \ref{totot} we have
$(\delta_4\delta_5)\beta_1^{-1}=\delta_1\delta_2$
and so 
$$((\delta_4\delta_5)^{-3}\delta_5(\delta_4\delta_5)^3)
\beta_1^{-1}=(\delta_1\delta_2)^{-4}\delta_2(\delta_1\delta_2)^4$$
On the other side, 
$$(\delta_1\delta_3)\beta_1^{-1}=\delta_1\delta_2
\delta_3\delta_5^{-1}\delta_4\delta_5\delta_3^{-1}\delta_5^{-1}\delta_4^{-1}
\delta_5\delta_4\delta_2^{-1}\delta_1^{-1}\delta_2$$
and using \ref{ftffff} and again \ref{totot}
$$(\delta_1\delta_3)\beta_1^{-1}=\delta_1\delta_2
\delta_3\delta_4\delta_5\delta_3^{-1}\delta_5^{-1}
\delta_2^{-1}\delta_1^{-1}\delta_2=\delta_1\delta_2
\delta_3\delta_4\delta_5\delta_3^{-1}
\delta_2^{-1}\delta_1^{-1}.$$
With the same method we have \begin{eqnarray*}
((\delta_1\delta_3)^{-1}\delta_3\delta_1\delta_3)\beta_1^{-1}&=&
\delta_1\delta_2\delta_3\delta_5^{-1}\delta_4^{-1}(\delta_3^{-1}\delta_5^{-1}
\delta_4^{-1}\delta_5\delta_4\delta_5\delta_3)\delta_4\delta_5\delta_3^{-1}
\delta_2^{-1}\delta_1^{-1}\\
&=&\delta_1\delta_2\delta_3\delta_5^{-1}\delta_4^{-1}
(\delta_5^{-1}\delta_4^{-1}
\delta_5\delta_4\delta_5)\delta_4\delta_5\delta_3^{-1}
\delta_2^{-1}\delta_1^{-1}\\
&=&\delta_1\delta_2\delta_3(\delta_4\delta_5)^{-2}
\delta_5(\delta_4\delta_5)^2\delta_3^{-1}
\delta_2^{-1}\delta_1^{-1}\end{eqnarray*}
So the relation is
$$\delta_3(\delta_4\delta_5)^{-2}\delta_5(\delta_4\delta_5)^2\delta_3^{-1}=
(\delta_1\delta_2)^{-5}\delta_2(\delta_1\delta_2)^5=
\delta_1\delta_2\delta_1^{-1}$$ where we have used 
\ref{otwtwo}.\\
Finally we have to write the relations 
$\Gamma\delta_i\Gamma^{-1}=
(\delta_i)\beta_1^2$, i.e.
\begin{eqnarray*}
\Gamma\delta_1\Gamma^{-1}&\!\!\!=\!\!\!&\Gamma\delta_2\Gamma^{-1}\delta_4^{-1}
\delta_2^{-1}\delta_1\delta_2\delta_4\Gamma\delta_2^{-1}\Gamma^{-1}
\hbox to 
0pt{\hspace*{13.6em}\begin{numera}\label{GoG}\end{numera}\hss}\\
\Gamma\delta_2\Gamma^{-1}&\!\!\!=\!\!\!&\delta_2^{-1}\delta_1\delta_2
\delta_3^{-1}\delta_5\delta_3\delta_2^{-1}\delta_1^{-1}\delta_2
\hbox to 
0pt{\hspace*{15.4em}\begin{numera}\label{GtwG}\end{numera}\hss}\\
\Gamma\delta_3\Gamma^{-1}&\!\!\!=\!\!\!&
\delta_4^{-1}\delta_2^{-1}\delta_1^{-1}\delta_2\delta_4\delta_2^{-1}
\delta_1\delta_2 
\delta_3^{-1}\delta_5^{-1}\delta_3\delta_5\delta_3
\delta_2^{-1}\delta_1^{-1}\delta_2\delta_4^{-1}\delta_2^{-1}\delta_1\delta_2
\delta_4
\ \ \begin{numera}\label{GthG}\end{numera}\\
\Gamma\delta_4\Gamma^{-1}&\!\!\!=\!\!\!&
\delta_4^{-1}\delta_2^{-1}\delta_1^{-1}\delta_2
\delta_4\delta_2^{-1}\delta_1\delta_2\delta_4
\hbox to 
0pt{\hspace*{15em}\begin{numera}\label{GfoG} 
\end{numera}\hss}\\
\Gamma\delta_5\Gamma^{-1}&\!\!\!=\!\!\!&\delta_5 
\hbox to
0pt{\hspace*{24.8em}
\begin{numera}\label{GfiG} \end{numera}\hss} 
\end{eqnarray*}
{\bf Step 4.}\\
We take as Shreier set for the left cosets of the 
kernel the set
$\{S_0=1,S_1=\delta_1\}$,
so applying the Reidemeister-Shreier method we get 
the generators
$\Delta=\delta_1^2$, $\Gamma$, 
$\sigma=\delta_1\Gamma\delta_1^{-1}$,
$A_i=\delta_1\delta_i$ and 
$B_i=\delta_i\delta_1^{-1}$ for $i\ugu2,3,4,5$.
For the relations we must rewrite the relations
\ref{fofififo},...,\ref{GfiG} and their conjugate by 
$\delta_1$ in terms of
the new generators. The rewriting process goes as 
follows (cf.
\cite{makaso}, pages 86-98): \begin{eqnarray*}
S_0\delta_1=S_1&&
S_1\delta_1=\Delta S_0\\
S_0\delta_i=B_iS_1&&
S_1\delta_i=A_iS_0\ \ \ \ {\rm for}\ i\ugu2,3,4,5\\
S_0\Gamma=\Gamma S_0&&
S_1\Gamma=\sigma S_1\end{eqnarray*}
We want to show that it suffices to rewrite only the 
relations
\ref{fofififo},...,\ref{GfiG}.\\
Observe that all our relations can be written in the 
form
$W\delta_iW^{-1}=\delta_k$ for a suitable word $W$. 
Assume that $\Gamma$
doesn't appear in the relation and do the rewriting 
after moding out by
the relations 
\begin{equaz}\De=B_iA_i=1.\label{canc1}\end{equaz}
Since $S_0\delta_i\ugu A_i^{-1}S_1$ and also 
$S_0\delta_i^{-1}\ugu
A_i^{-1}S_1$, if we write
$\DS 
W=\prod_{\lambda=1}^h\delta_{j_{\lambda}}^{\pm1}$, 
the rewriting of
$W\delta_iW^{-1}\delta_k^{-1}$ is given by
$$A_{j_1}^{-1}A_{j_2}\cdots A_i^{\pm1}\cdots 
A_{j_1}^{-1}A_k$$
(note that $A_1\ugu1$). The rewriting of the same 
relation conjugated by
$\delta_1$ yields instead
$$A_{j_1}A_{j_2}^{-1}\cdots A_i^{\mp1}\cdots 
A_{j_1}A_k^{-1}.$$
We get thus two relations of respective form 
$UA_k=1$, $U^{-1}A_k^{-1}=1$,
which are obviously equivalent.\\
If instead $\Gamma$ appears in the relation, we have 
one of the
\ref{GoG},...,\ref{GfiG} which are of the form
$\Gamma\delta_i\Gamma^{-1}=W\delta_iW^{-1}$ where we 
can in fact assume that
$\Gamma$ doesn't appear in the word $W$.\\ The 
rewriting of
$\Gamma\delta_i\Gamma^{-1}W\delta_i^{-1}W^{-1}$ 
yields, again a relation of
the form
$$\Gamma A_i^{-1}\sigma^{-1}U^{-1}=1,$$ whereas the 
rewriting of the
conjugate by $\delta_1$ gives a relation $$\sigma 
A_i\Gamma^{-1}U=1,$$ which
is an equivalent relation.\\
For convenience of notation we shall keep the 
generators $B_i=A_i^{-1}$.\\
To calculate $\pi_1(Z_2^{\#})$ we must add the 
rewriting of
$(\prod_{i=1}^5\delta_i)^2=1$ which gives
$$A_2B_3A_4B_5B_2A_3B_4A_5=1.$$
We have thus that $\pi_1(Z_2^{\#})$ is generated by 
$A_2$, $A_3$, $A_4$,
$A_5$, $\Gamma$ and $\sigma$ and has the
following set of relations
\begin{equaz}\label{3)'} (B_4A_5)^6=(B_5A_4)^3 
\end{equaz}
\begin{equaz}\label{1)'} B_3A_2=B_2A_3\end{equaz}
\begin{equaz}\label{4)'} B_5=B_2^3 \end{equaz}
\begin{equaz}\label{2)'} B_3^3=(B_5A_4)^6B_5 \end{equaz}
\begin{equaz}\label{5)'} A_2^{12}=1 \end{equaz}
\begin{equaz}\label{6)'} 
B_5A_3B_5A_4B_5A_4=B_4A_5B_4A_5B_3A_5 \end{equaz}
\begin{equaz}\label{7)'} 
B_5=B_2^2A_4B_5A_3B_5A_4B_2^2 \end{equaz}
\begin{equaz}\label{8)'} B_3B_2B_3=(B_5A_4)^4B_5 
\end{equaz}
\begin{equaz}\label{12)'} \sigma 
A_2^2\Gamma^{-1}=A_4B_2^2A_4 \end{equaz}
\begin{equaz}\label{10)'}\Gamma 
B_2\sigma^{-1}=B_2^2A_3B_5A_3B_2^2\end{equaz}
\begin{equaz}\label{13)'}\Gamma B_3\sigma^{-1}=
B_4A_2^2B_4A_2^2B_3A_5B_3A_5B_3A_2^2B_4A_2^2B_4\end{equaz}
\begin{equaz}\label{11)'}\Gamma 
B_4\sigma^{-1}=B_4A_2^2B_4A_2^2B_4\end{equaz}
\begin{equaz}\label{9)'} \Gamma 
B_5\sigma^{-1}=B_5\end{equaz}
\begin{equaz}\label{14)'} 
A_2B_3A_4B_5B_2A_3B_4A_5=1\end{equaz}
where $B_i=A_i^{-1}$.\\
Let's reduce this presentation. Using \ref{4)'} 
relation \ref{7)'}
becomes 
\begin{equaz}B_4A_2B_4=B_5A_3B_5\label{7)''}\end{equaz} and with
this \ref{6)'}
becomes $$B_2^4=1$$ which implies \ref{5)'}, and 
changes \ref{4)'} into
$$B_5=A_2.$$
Moreover, using \ref{8)'} and the last equation, 
relation
\ref{2)'} gives $$(A_2A_4)^2=A_3A_2B_3^2$$ and with 
this, using also
\ref{1)'}, \ref{8)'} becomes
\begin{equaz}\label{8)''}B_3B_2=A_2A_3\end{equaz}
thus transforming \ref{7)''} into
$$B_3=B_4A_2B_4$$ which allows us to delete the 
generator $A_3$.
Upon substituting the expressions of $A_5$ and $A_3$ 
into \ref{14)'} and
\ref{8)''} we have $$A_2A_4=A_4A_2$$ $$A_4^4=1.$$
We can then see that the relations 
\ref{1)'},...,\ref{14)'}
are equivalent to the following
$$A_5=A_2^{-1}\ \ \ \ \ \ \
A_3=A_2^{-1}A_4^2$$
$$A_2^4=A_4^4=1\ \ \ \ \ \ \
A_2A_4=A_4A_2$$
$$\sigma A_2^2\Gamma^{-1}=A_2^2A_4^2$$
$$\Gamma A_2^{-1}=A_2^{-1}\sigma$$
$$\Gamma A_2A_4^2=A_2A_4^2\sigma$$
$$\Gamma A_4^{-1}=A_4\sigma$$
$$\Gamma A_2=A_2\sigma.$$
{\bf Step 5.}\\
We take as Shreier set for the left cosets of the 
kernel the set
$$\{R_i\ugu\Gamma^i\ |\ i\ugu0,1\vir m\men1\}$$ and 
we apply the
Reidemeister-Shreier method.\\
The generators are $\hat\Gamma=\Gamma^m$,
$A_{2,i}=\Gamma^iA_2\Gamma^{-i}$, 
$A_{4,i}=\Gamma^iA_4\Gamma^{-i}$ for
$i\ugu0\vir m\men1$ 
$\sigma_i=\Gamma^i\sigma\Gamma^{-(i+1)}$ for
$i\ugu0\vir m\men2$ and 
$\sigma_{m-1}=\Gamma^{m-1}\sigma$.\\
For the rewriting process we have
\begin{eqnarray*}R_iA_j&=&A_{j,i}R_i\ \ \ \ {\rm 
for}\ j\ugu2,4\ i\ugu0\vir
m-1\\
R_i\Gamma&=&R_{i+1}\ \ \ \ \ {\rm for}\ i\ugu0\vir m-2\\
R_{m-1}\Gamma&=&\hat\Gamma R_0\\
R_i\sigma&=&\sigma_i R_{i+1}\ \ \ \ {\rm for}\ 
i\ugu0\vir m-2\\
R_{m-1}\sigma&=&\sigma_{m-1}R_0.
\end{eqnarray*}
Thus, taking indices $i$ (mod$m$) and adding (as we 
must) the relation
$\hat\Gamma=1$, we obtain the relations
$$A_{2,i}^4=A_{4,i}^4=1$$
$$A_{2,i}A_{4,i}=A_{4,i}A_{2,i}$$
$$\sigma_iA_{2,i+1}^2=A_{2,i}^2A_{4,i}^2$$
$$A_{2,i+1}^{-1}=A_{2,i}^{-1}\sigma_i$$
$$A_{2,i+1}A_{4,i+1}^2=A_{2,i}A_{4,i}^2\sigma_i$$
$$A_{4,i+1}^{-1}=A_{4,i}\sigma_i$$
$$A_{2,i+1}=A_{2,i}\sigma_i.$$
To simplify this presentation we write
\begin{eqnarray*}\sigma_i&=&A_{2,i}^2A_{4,i}^2A_{2,i+1}^2\\
&=&A_{2,i}A_{2,i+1}^{-1}\\
&=&A_{4,i}^2A_{2,i}^{-1}A_{2,i+1}A_{4,i+1}^2=
A_{2,i}^{-1}A_{4,i}^2A_{4,i+1}^2A_{2,i+1}\\
&=&A_{4,i}^{-1}A_{4,i+1}^{-1}\\
&=&A_{2,i}^{-1}A_{2,i+1}\end{eqnarray*}
 From the last and the second equations we get
$$A_{2,i}^2=A_{2,i+1}^2=A_{2,0}^2$$ and from the 
first one,
remembering that $A_{2,i}$ commutes with $A_{4,i}$ 
and that $A_{2,0}^4=1$,
$$\sigma_i=A_{4,i}^2.$$
The fourth equation then
gives $$A_{4,i}=A_{4,i+1}=A_{4,0}$$ which makes the 
last and the
third relations equivalent. These two cancellation 
relations
enable us to delete all the generators
$\sigma_j$ and $A_{4,i}$ for $i\ugu1\vir m-1$.\\
We may rewrite the five relations above as
\begin{eqnarray*}\sigma_i&=&A_{4,0}^2\\
A_{4,i}&=&A_{4,0}\\
A_{4,0}^2&=&A_{2,i}^{-1}A_{2,i+1}\\
A_{4,0}^2&=&A_{2,i}A_{2,i+1}^{-1}.
\end{eqnarray*}
Clearly the last two equations are equivalent and give
\begin{equaz}A_{2,2i}=A_{2,0}\ \ \ \ \ \
A_{2,2i+1}=A_{2,0}A_{4,0}^2\label{abv}.\end{equaz}
Moreover, if we add the relation $\sigma^m=1$, that 
in the generators of
$\pi_1(X_{2m}^{\#})$ reads out as 
$\sigma_0\sigma_1\cdots\sigma_{m-1}=1$, we
get $A_{4,0}^{2m}=1,$ i.e., if $m$ is odd, 
$A_{4,0}^2=1$, while if $m$ is
even we have no new relations. Observe that this is 
in accordance with the
fact that in \ref{abv} the index is cyclic mod$(m)$.\\
Summing up, we have a commutative group with only 
two generators, namely
$a=A_{2,0}$ and $b=A_{4,0}$, such that $a^4=1$ and 
$b^4=1$ if $m$ is even,
$b^2=1$ if $m$ is odd, i.e.
$$\pi_1(X^{\#}_{2k+2})=\zeta_4\times\zeta_4$$ if $k$ 
is odd and
$$\pi_1(X^{\#}_{2k+2})=\zeta_4\times\zeta_2$$ if $k$ 
is even.

\end{document}